\documentclass[12pt]{msml2021} 

\title[Optimal Policies for a Pandemic: Modeling and Algorithm]{Optimal Policies for a Pandemic: A Stochastic Game Approach and a Deep Learning Algorithm}
\usepackage{times}
\usepackage{bm,amsmath}
\usepackage{mathtools,empheq,algorithm,algorithmic}
\usepackage{physics,comment}
\usepackage{graphics}
\mathtoolsset{showonlyrefs}

\DeclareMathOperator*{\argmin}{arg\,min} 
 
\newcommand{\algorithmicdoinparallel}{\textbf{do in parallel}}
\makeatletter
\AtBeginEnvironment{algorithmic}{%
  \newcommand{\FORALLP}[2][default]{\ALC@it\algorithmicforall\ #2\ %
    \algorithmicdoinparallel\ALC@com{#1}\begin{ALC@for}}%
}

\newcommand{\cN}{\mathcal{N}}

\newcommand{\bX}{\bm{X}}
\newcommand{\bx}{{\bm{x}}}

\newcommand{\bW}{\bm{W}}
\newcommand{\bl}{\bm{\ell}}
\newcommand{\bh}{\bm{h}}

\newcommand{\bs}{\bm{s}}
\newcommand{\be}{\bm{e}}
\newcommand{\bi}{\bm{i}}
\newcommand{\tj}{{j'}}
\newcommand{\balpha}{\bm{\alpha}}

\newcommand{\mc}[1]{\mathcal{#1}}

\newcommand{\MCO}{\mathcal{O}}

\newcommand{\EE}{\mathbb{E}}

\newcommand{\RR}{\mathbb{R}}

\newcommand{\ltwonorm}[1]{\left\lVert#1\right\rVert_2}

\newcommand{\dm}[4]{{\beta^{{#1}{#2}}s_{#1}i_{#2}\left(\pdv{V^{{#3}, {#4}}}{e_{#1}} - \pdv{V^{{#3}, {#4}}}{s_{#1}}\right)}}
\newcommand{\transpose}{^{\operatorname{T}}}

\newcommand{\half}{\frac{1}{2}}

\newcommand{\eps}{\epsilon}

\newcommand{\ud}{\,\mathrm{d}}

\newtheorem{theo}{Theorem}[section]

\newtheorem{defi}[theo]{Definition}

\msmlauthor{%
 \Name{Yao Xuan} \Email{yxuan@math.ucsb.edu}\\
 \Name{Robert Balkin} \Email{rbalkin@ucsb.edu}\\
 \addr Department of Mathematics, University of California, Santa Barbara, CA 93106-3080, USA 
 \AND
 \Name{Jiequn Han} \Email{jiequnh@princeton.edu}\\
 \addr Department of Mathematics, Princeton University, Princeton, NJ 08544-1000, USA 
 \AND
  \Name{Ruimeng Hu} \Email{rhu@ucsb.edu}\\
 \addr Department of Mathematics and Department of Statistics and Applied Probability, University of California, Santa Barbara, CA 93106-3080, USA
\AND
  \Name{Hector D. Ceniceros} \Email{ceniceros@ucsb.edu}\\
 \addr Department of Mathematics, University of California, Santa Barbara, CA 93106-3080, USA 
}

\makeatletter
\let\Ginclude@graphics\@org@Ginclude@graphics
\makeatother

\begin{document}

\maketitle

\begin{abstract}%
Game theory has been an effective tool in the control of disease spread and in suggesting optimal policies at both individual and area levels. In this paper, we propose a multi-region SEIR model based on stochastic differential game theory, aiming to formulate  optimal regional policies for infectious diseases. Specifically, we enhance the standard epidemic SEIR model by taking into account the social and health policies issued by multiple region planners. This enhancement makes the model more realistic and powerful. However, it also introduces a formidable computational challenge due to the high dimensionality of the solution space brought by the presence of multiple regions. This significant numerical difficulty of the model structure motivates us to generalize the deep fictitious algorithm introduced in [Han and Hu, MSML2020, pp.221--245, PMLR, 2020] and develop an improved algorithm to overcome the curse of dimensionality. We apply the proposed model and algorithm to study the COVID-19 pandemic in three states: New York, New Jersey and Pennsylvania. The  model parameters are estimated from real data posted by the Centers for Disease Control and Prevention (CDC). We are able to show the effects of the lockdown/travel ban policy on the spread of COVID-19 for each state and how their policies affect each other.

\end{abstract}


\section{Introduction}

The pandemic of coronavirus disease 2019 (COVID-19) has brought a huge impact on our lives. Based on the CDC Data Tracker, as of early December 2020, there have been more than 15 million confirmed cases of infection and more than 290 thousand cases of death in the United States. Needless to say,  the economic impact has also been catastrophic, resulting in unprecedented unemployment and the bankruptcy of many restaurants, recreation centers, shopping malls, etc.

In a classic, compartmental epidemiological model each individual is assigned a label, {\it e.g.}, \textbf{S}usceptible, \textbf{E}xposed, \textbf{I}nfectious, \textbf{R}emoved, \textbf{V}accinated. The labels' order shows the flow patterns between the compartments (SIR, SEIR, SIRV models). Other approaches include network models, which explicitly include the interaction of individuals, in addition to the modeling of each individual's dynamics, and agent-based models that are useful in informing decision making when accurately calibrated.
Moreover, the consideration of pharmaceutical and/or non-pharmaceutical intervention policies naturally couples game theory to epidemiological models by controlling when and how the game is played in such models. For example, in some early studies \cite{bauch2003group,bauch2004vaccination}, one can use non-repeated games to incorporate  game theory into modeling at the individual level, where individuals (known as ``players'' in the game theory) maximize their gain by weighing the costs and benefits of different strategies. We refer to the review paper \cite{gamereviewpandemic} and the references therein for more details.

Differential games, initiated by \cite{Isaacs1965}, as an offspring of game theory and optimal control, provide modeling and analysis of conflict in the context of a dynamical system. They have been intensively employed across many disciplines, including management science,  economics, social science, biology, military, etc. One of the core objectives in differential games is to compute Nash equilibria that refer to strategies by which no player has an incentive to deviate. However, a major {bottleneck} comes from the notorious intractability of $N$-player games, and the direct computation of Nash equilibria is extremely time-consuming and memory demanding. In a series of recent works by \cite{Hu2:19,HaHu:19,han2020convergence,hanhulong:21}, the deep fictitious play (DFP) theory and algorithms were developed for stochastic differential games (SDG) with a large number of heterogeneous players. The DFP framework embeds the fictitious play idea, introduced by \cite{Br:49,Br:51}, into designed architectures of deep neural networks to produce accurate and parallelizable algorithms with convergence analysis, and resolve the intractability issue (curse of dimensionality) caused by the complex modeling and underlying high-dimensional space in SDG.   

Building from the DFP theory and algorithms for computing Nash equilibria in SDG, we propose here  to strengthen the classical SEIR model by taking into account the social and health policies issued by multiple region planners. We call this new model
a stochastic multi-region SEIR model because it couples the stochastic differential game theory with the SEIR model, making it  more realistic and powerful. The computational challenge introduced by the high-dimensionality of the multi-region solution space is addressed by generalizing the deep fictitious algorithm proposed by \cite{HaHu:19}. This new approach leads to an enhanced deep fictitious play algorithm to overcome the curse of dimensionality and further reduce the computational complexity. To showcase the performance of the proposed model and algorithm, we apply them to a case study of the COVID-19 pandemic in three states: New York (NY), New Jersey (NJ), and Pennsylvania (PA). We present the optimal lockdown policy corresponding to the Nash equilibrium of the multi-region SEIR model. We remark that our work is not to predict a pandemic, but to provide a game-themed framework, a deep learning algorithm and  possible outcomes for competitive region planners. We hope that information can provide some qualitative guidance for policymakers on the impact of certain policies. The parameters used in the numerical experiments are based on the current knowledge of the Coronavirus, which is still under development. In practice, at the beginning of the Coronavirus, a governor might not be able to trace the infections fully, and infected cases may not be fully identified. All these may lead to the inaccuracy of parameter estimations despite our using of the best available data. Therefore, it may be hard to match the predicted results with practical observation.

The rest of the paper is organized as follows. In Section~\ref{sec:model}, we propose a novel multi-region epidemiological model and explain how social and health policies issued by region planners are integrated into  dynamics, leading to a game feature of the problem. We explain the numerical challenge in Section~\ref{sec:algorithm} and propose an enhanced deep fictitious play algorithm for memory and computational efficiency.  Section~\ref{sec:case} focuses on a NY-NJ-PA case study with detailed discussions on parameter choices and the resulted optimal policies. We present some concluding remarks in Section~\ref{sec:conclusion} and provide technical details in the appendices.

\section{Mathematical modeling: A multi-region SEIR model}\label{sec:model}
We consider a pandemic spreading in $N$ geographical regions, and each planner controls the loss of her region by implementing some policies. We aim to study how the region planners' policies affect each other, and the equilibrium policies. 

Let us start with a modified version of the very-known epidemic SEIR model (cf.  \cite{dynamicalliu}), where each region's population is assigned to compartments with four labels: \textbf{S}usceptible, \textbf{E}xposed, \textbf{I}nfectious, and \textbf{R}emoved. Individuals with different labels denote $S$: those who are not yet infected; $E$: who have been infected but are not yet infectious themselves; $I$: who have been infected and are capable of spreading the disease to those in the susceptible category, and $R$: who have been infected and then removed from the disease due to recovery or death. The region planners can issue certain policies to mitigate the pandemic, for instance, policies that can help reduce the transmission rates and death rates. Mathematically, denote by $S^n_t, E^n_t, I^n_t, R^n_t$ the \emph{proportion} of population in the four compartments of region $n$ at time $t$. We consider the following stochastic multi-region SEIR model: 
\begin{align}
&\ud S^n_t = -\sum_{k = 1}^N \beta^{nk} S^n_t  I^k_t  (1-\theta \ell^n_t)(1-\theta \ell^k_t) \ud t - v(h^n_t)S^n_t \ud t - \sigma_{s_n} S^n_t \ud W^{s_n}_t,\label{def:St}\\
&\ud E^n_t = \sum_{k = 1}^N \beta^{nk} S^n_t  I^k_t  (1-\theta \ell^n_t)(1-\theta \ell^k_t) \ud t - \gamma E^n_t \ud t + \sigma_{s_n} S^n_t \ud W^{s_n}_t - \sigma_{e_n} E^n_t \ud W^{e_n}_t,\label{def:Et} \\
&\ud I^n_t  = (\gamma E^n_t - \lambda(h^n_t)I^n_t) \ud t + \sigma_{e_n} E^n_t \ud W^{e_n}_t, \label{def:It} \\
&\ud R^n_t  = \lambda(h^n_t)I^n_t \ud t + v(h^n_t) S^n_t \ud t, \quad n \in \mc{N} := \{1, 2, \ldots, N\}, \label{def:Rt}
\end{align}
where $\bm\ell_t \equiv (\ell^1_t, \ldots, \ell^N_t)$ and $\bm h_t \equiv (h^1_t, \ldots, h^N_t)$ are policies chosen by the region planners at time $t$.  Each planner $n$ seeks to minimize its region's cost within a period $[0,T]$:
\begin{equation}\label{def:J}
J^n(\bm \ell, \bm h) := \mathbb{E} \bigg[ \int_0^T e^{-rt} P^n\big[ (S^n_t + E^n_t + I^n_t)  \ell^n_t w + a(\kappa I^n_t \chi + p I^n_t c)\big]  + e^{-rt}\eta (h^n_t)^2\ud t \bigg].
\end{equation}

We now give detailed description of this model \eqref{def:St}--\eqref{def:J}:
\begin{enumerate}
	\item[S:] $\beta^{nk}$ denotes the average number of contacts per person per time. The transition rate between $S^n$ and $E^n$ due to contacting infectious people in the region $k$ is proportional to the fraction of those contacts between an infectious and a susceptible individual, which result in the susceptible one becoming infected, {\it i.e.}, $\beta S^n_tI^k_t$. Although some regions may not be geographically connected, the transmission between the two is still possible due to air travels but is less intensive than the transmission within the region,  {\it i.e.}, $\beta^{nk}>0$ and $\beta^{nn} \gg \beta^{nk}$ for all $k \neq n$. 
	
	$\ell^n_t \in [0, 1]$ denotes the decision of the planner $n$ on the fraction of population being locked down at time $t$. We assume that those in lockdown cannot be infected. However, the policy may only be partially effective as essential activities (food production and distribution, health, and basic services) have to continue. Here we use $\theta \in [0,1]$ to measure this effectiveness,  and the transition rate under the policy $\bm \ell$ thus become $\beta^{nk} S^n_t  I^k_t  (1-\theta \ell^n_t)(1-\theta \ell^k_t)$. The case $\theta = 1$ means the policy is fully effective. 
	
	$h^n_t \in [0,1]$ denotes the effort the planner $n$ decide to put into the health system, which we refer as \emph{health policy}. It will influence the vaccination availability $v(\cdot)$ and the recovery rate $\lambda(\cdot)$ of this model.
	
	$v(h^n_t)$ denotes the vaccination availability of region $n$ at time $t$. Once vaccinated, the susceptible individuals $v(h^n_t)S^n_t$ become immune to the disease, and join the removed category $R^n_t$. We model it as an increasing function of $h^n_t$, and if the vaccine has not been developed yet, we can define $v(x) = 0$ for $x \leq \overline h$.

	\item[E:] $\gamma$ describes the latent period when the person has been infected but not infectious yet. It is the inverse of the average latent time, and we assume $\gamma$ to be identical across all regions. The transition between $E^n$ and $I^n$ is proportional to the fraction of exposed, {\it i.e.}, $\gamma E^n_t$. 
	
	\item[I:]   $\lambda(\cdot)$ represents the recovery rate. For the infected individuals, a fraction $\lambda(h^n)I^n$ (including both death and recovery from the infection) joins the removed category $R^n$ per time unit.
	The rate is determined by the average duration of infection $D$. We model the duration (so does the recovery rate) related to the health policy $h^n_t$ decided by its planner. The more effort put into the region ({\it i.e.}, expanding hospital capacity, creating more drive-thru testing sites), the more clinical resources the region will have and the more resources will be accessible by patients, which could accelerate the recovery and slow down death. The death rate, denoted by $\kappa(\cdot)$, is crucial for computing the cost of the region $n$; see the next item.
	 
	\item[Cost:] Each region planner faces four types of cost. One is the economic activity loss due to the lockdown policy, where $w$ is the productivity rate per individual, and $P^n$ is the population of the region $n$. The second one is due to the death of infected individuals. Here $\kappa$ is the death rate which we assume for simplicity to be constant, and $\chi$ denotes the cost of each death.  The hyperparameter $a$ describes how planners weigh deaths and infections comparing to other costs. 
	The third one is the inpatient cost, where $p$ is the hospitalization rate, and $c$ is the cost per inpatient day. The last term $\eta (h^n_t)^2$ is the grants putting into the health system. We choose a quadratic form to account for diminishing marginal utility (view it from $\eta (h^n_t)^2$ to $h^n_t$). All costs are discounted by an exponential function $e^{-rt}$, where $r$ is the risk-free interest rate, to take into account the time preference.
	Note that region $n$'s cost depends on all regions' policies $(\bm \ell, \bm h)$, as $\{I^k, k \neq n\}$ appearing in the dynamics of $S^n$. Thus we write is $J^n(\bm \ell, \bm h)$.

\end{enumerate}
The choices of epidemiological parameters will be discussed in Section~\ref{sec:parameter}. Next, we summarize the key assumptions in the above model:
\begin{enumerate}
    \item The dynamics of an epidemic are much faster than the vital (birth and death) dynamics. So vital dynamics are omitted in the above model.
	\item The planning is of a short horizon and will be adjusted frequently as the epidemic develop. For simplicity, we assume this is no migration between regions over the time $[0,T]$. 
	\item Individuals who once recovered from the disease, are immune and free of lockdown policy. 
	\item The dynamics obeys the conservation law: $S^n_t + E^n_t + I^n_t + R^n_t = P^n$. This means that the process $R^n$ is redundant.
    \item The dynamics of $S$, $E$ and $I$ are subjected to random noise, to account for the noise introduced during data recording, false-positive/negative test results, exceptional cases when recovered individuals become susceptible again, minor individual differences in the latent period, etc.
    \item Individuals who are not under lockdown have the same productivity, no matter their categories. We assume this for simplicity remark that this can be improved by assigning different productivity to individuals with or without symptoms.

\end{enumerate}

The above modeling and objectives can be viewed as a stochastic differential game between $N$ players\footnote{Henceforth, we shall use \emph{planner} and \emph{player} interchangeably}. 
Here, we view the whole problem as a non-cooperative game, as many regions make decisions individually and indeed even compete for scarce resources (frontline workers, personal protective equipment, etc.) during the outbreak. 
Each player $n$ controls her states $(S^n, E^n, I^n, R^n)$ through her strategy $(\ell^n, h^n)$ in order to minimize the associated cost $J^n$. The optimizers then are interpreted as the optimal lockdown policy and optimal effort putting into the health system. 

For a non-cooperative game, one usually refers to Nash equilibrium as a notion of optimality. For completeness, we review the definition here.
\begin{defi}
A Nash equilibrium is a tuple $(\bm \ell^\ast, \bm h^\ast) = (\ell^{1, \ast}, h^{1, \ast}, \ldots, \ell^{N, \ast}, h^{N, \ast} ) \in \mathbb{A}^N$ such that 
\begin{equation}
\forall n \in \mc{N}, \emph{ and } (\ell^n, h^n) \in \mathbb{A}, \quad J^n(\bm \ell^\ast, \bm h^\ast) \leq J^n((\bm \ell^{-n, \ast}, \ell^n), (\bm h^{-n, \ast}, h^n)),
\end{equation}
where $\bm \ell^{-n,\ast}$ represents strategies of players other than the $n$-th one:
\begin{equation}
    \bm\ell^{-n, \ast} := [\ell^{1, \ast}, \ldots, \ell^{n-1, \ast}, \ell^{n+1, \ast}, \ldots, \ell^{N,\ast}] \in \mathbb{A}^{N-1},
\end{equation}
$\mathbb{A}$ denotes the set of admissible strategies for each player and  $\mathbb{A}^N$ is the produce of $N$ copies of $\mathbb{A}$. For simplicity, we have assumed all players taking actions in the same space.
\end{defi}

In the sequel, to fix the notations, we shall use
\begin{itemize}
    \item  a regular character with a superscript $n$ for an object from player $n$;
    \item a boldface character for a collection of objects from all players, {\it i.e.}, $\bm S_t \equiv [S_t^1, \ldots, S_t^N]\transpose$; 
    \item a boldface character with a superscript $-n$ for a collection of objects from all players except $n$, {\it i.e.}, $\bm S_t^{-n} \equiv [S_t^1, \ldots, S_t^{n-1}, S_t^{n+1}, \ldots,  S_t^N]\transpose$.
\end{itemize}

A Markovian Nash equilibrium is a Nash equilibrium defined above with $\mathbb{A}$ being the set of Borel measurable functions: $(\ell, h): [0, T] \times \RR^{3N} \to [0,1]^2$. In other words, the policies $(\ell_t^n, h_t^n)$ at time $t$ are functions of the time $t$ and the current values of all players' state processes $(\bm S_t, \bm E_t, \bm I_t)$. We omit the dependence on $\bm R_t$ as it is redundant as a consequence of the conservation law.

We derive below  the Hamilton-Jacobi-Bellman (HJB) equations characterizing the Markovian Nash equilibrium. To simplify the notation,  we first rewrite the dynamics of $(\bm S_t, \bm E_t, \bm I_t)$ defined in \eqref{def:St}--\eqref{def:It} into a vector form $\bX_t \equiv [\bm S_t,\bm E_t, \bm I_t]\transpose \equiv [S_t^1, \cdots, S_t^N, E_t^1, \cdots, E_t^N, I_t^1, \cdots, I_t^N]\transpose \in \RR^{3N}$. Again, we shall drop the redundant process $\bm R_t$. The dynamics of $\bX_t$ reads:
\begin{equation}\label{def:Xt}
    \ud \bX_t = b(t, \bX_t, \bm \ell(t, \bX_t), \bm h(t, \bX_t)) \ud t  + \Sigma(\bX_t)\ud \bm W_t, 
\end{equation}
where $b, \Sigma$ are deterministic functions in $\RR^{3N}$ and  $\RR^{3N \times 2N}$, and $\{\bm W_t\}_{0\leq t \leq T}$ is a $2N$-dimensional standard Brownian motion. Each player $n$ aims to minimize the expected running cost
\begin{equation}\label{def:cost}
 \EE\left[\int_0^T f^n(t, \bX_t, \ell^n(t, \bX_t), h^n(t, \bX_t)) \ud t \right].  
\end{equation}
We defer the specific definitions of $b, \Sigma, \bm W$, and $f^n$ to Appendix~\ref{app:details-sec2} to facilitate the exposition. We now define the value function of player $n$ by
\begin{equation}
    V^n(t, \bx) = \inf_{(\ell^n, h^n) \in \mathbb{A}} \EE\left[\int_t^T f^n(s, \bX_s, \ell^n(s, \bX_s), h^n(s, \bX_s))\ud s\vert \bX_t = \bx\right].
\end{equation}
By dynamic programming, it solves the following HJB system
\begin{equation} \label{def:HJB}
\begin{dcases}
\partial_t V^n + \inf_{(\ell^n, h^n) \in [0,1]^2} H^n(t,\bm x,(\bl, \bh)(t, \bx), \nabla_{\bx} V^n)  + \half \text{Tr}(\Sigma(\bx)\transpose \text{Hess}_{\bx} V^n \Sigma(\bx)) = 0,\\
V^n(T,\bx) = 0, \quad n \in \mc{N},
\end{dcases}
\end{equation}
where $H^n$ is the usual Hamiltonian defined by 
\begin{equation}\label{def:H}
    H^n(t,\bm x, \bl, \bh, \bm p) =  b(t,\bm x, \bl, \bh) \cdot \bm p + f^n(t, \bm x, \ell^n, h^n),
\end{equation}
$\partial_t$ denotes the time derivative, $\nabla_{\bx} V$ and~$\text{Hess}_{\bx} V$ denote
the gradient and the Hessian of the function $V$ with respect to $\bm x$, respectively,  and $\text{Tr}$ stands for the trace of a matrix. 

Finding the optimal policies for $N$ regions is equivalent to solving $N$-coupled $3N+1$ dimensional nonlinear equations \eqref{def:HJB}. For example, when $N = 3$, each PDE is 10-dimensional and conventional methods start to lose their efficiency. The recently proposed deep learning algorithm in \cite{HaHu:19}, known as deep fictitious play (DFP), has shown excellent numerical performance in solving high-dimensional,  coupled HJB equations with convergence analysis \citep{han2020convergence}. In the next section, we will first review and then propose an enhanced version of DFP to tackle some new issues.

\section{Numerical methodology: Enhanced deep fictitious algorithm}\label{sec:algorithm}
We first briefly review the deep fictitious play (DFP) algorithm for solving equation \eqref{def:HJB} and we refer 
readers to \cite{HaHu:19} for full details. With the idea of fictitious play, DFP recasts the $N$-player game into $N$ decoupled optimization problems, which are solved repeatedly stage by stage. Each individual problem is solved by the deep BSDE method \citep{EHaJe:17,HaJeE:18}. The algorithm starts with some initial guess $(\bl^0, \bh^0)$, where the superscript $0$ stands for stage $0$. At the $(m+1)^{th}$ stage, given the optimal policies $(\bl^m, \bh^m)$ at the previous stage, the algorithm solves the following PDEs
\begin{equation} \label{def:HJB-DFP}
\begin{dcases}
\partial_t V^{n, m+1} + \inf_{(\ell^n, h^n) \in [0,1]^2} H^n(t,\bm x,(\ell^n, \bl^{-n,m}, h^n, \bh^{-n,m})(t, \bx), \nabla_{\bx} V^{n, m+1})  \\
\hspace{50pt}+ \half \text{Tr}(\Sigma(\bx)\transpose \text{Hess}_{\bx} V^{n, m+1} \Sigma(\bx)) = 0,\\
V^{n, m+1}(T,\bx) = 0, \quad n \in \mc{N},
\end{dcases}
\end{equation}
and obtains the $(m+1)^{th}$ stage's optimal strategy by:
\begin{equation}\label{def:alphaast}
(\ell^{n, m+1}, h^{n, m+1})(t, \bm x) = \argmin_{(\ell^n, h^n) \in [0,1]^2} H^n(t, \bm x, (\ell^n, \bl^{-n,m}, h^n, \bh^{-n,m})(t, \bx), \nabla_{\bx} V^{n, m+1}(t, \bx)).
\end{equation}
Here,  $(\bl^{-n,m}, \bh^{-n,m})$ stands for others' optimal policies from the $m^{th}$ stage and are considered to be fixed functions when solving the PDE at the current stage. In the sequel, to simplify notations we omit the stage label $m$ in the superscript when there is no risk of confusion. To solve \eqref{def:HJB-DFP} at each stage, it is first rewritten in the DFP as
\begin{align}\label{def:PDE-FK}
\partial_t V^n + \half \text{Tr}(\Sigma(\bx)\transpose \text{Hess}_{\bx} V^n \Sigma(\bx)) &+ \mu^n(t, \bx; \bl^{-n}, \bh^{-n})\cdot \nabla_{\bx} V^n \\
&+g^n(t, \bx, \Sigma(\bx)\transpose\nabla_{\bx} V^n; \bl^{-n}, \bh^{-n})=0,
\end{align}
 with some functions $\mu^n$ and $g^n$. The solution is then approximated by solving the equivalent BSDE $(\bX_t^n, Y_t^n, Z_t^n)\in \RR^{3N} \times \RR \times \RR^{2N}$: 
\begin{empheq}[left=\empheqlbrace]{align}
    &\bX_t^n = \bx_0 + \int_{0}^{t} \mu^n(s, \bX_s^n; (\bl^{-n}, \bh^{-n})(s, \bX_s^n)) \ud s + \int_{0}^{t}\Sigma(\bX_s^n)\, \mathrm{d}\bW_s, \label{eq:BSDE_forward} \\
    &Y_t^n =  \int_{t}^{T}g^n(s, \bX_s^n, Z_s^n; (\bl^{-n}, \bh^{-n})(s, \bX_s^n))\ud s - \int_{t}^{T}(Z_s^n)\transpose\, \mathrm{d}\bW_s, \label{eq:BSDE_backward}
\end{empheq}
in the sense of (cf. \cite{PaPe:92,ElPeQu:97,PaTa:99})
\begin{equation}\label{eq:FK relation}
    Y_t^n = V^n(t, \bX_t^n) \quad\text{and}\quad Z_t^n = \Sigma(\bX_t^n)\transpose\nabla_{\bx} V^n(t, \bX_t^n).
\end{equation}
The high-dimensional BSDE \eqref{eq:BSDE_forward}--\eqref{eq:BSDE_backward} is tackled by the deep BSDE method proposed in \cite{EHaJe:17,HaJeE:18}. 

In \cite{HaHu:19}, the algorithm solves the BSDE by parametrizing $V^n(t, \bx)$ using neural networks (NN) and then obtains the approximate optimal policy by plugging the NN outputs into \eqref{def:alphaast}. For memory efficiency, the algorithm only stores the NNs' parameters at the current and the one-step previous stages. This strategy works well for games like the linear-quadratic game, but it would be ineffective  if  $\bl^{-n}$ or $\bh^{-n}$ explicitly appears in the minimizer in  \eqref{def:alphaast}. In this case, when evaluating others' strategy $\bl^{-n}$ or $\bh^{-n}$ at stage $m$, it does not only need NNs at stage $m$ but also at stages $m-1$, $m-2$, ..., $0$. This means one needs to store NNs' parameters for all the previous stages from $1, \ldots, m$, and evaluate the associated output.
Therefore, the time complexity of evaluating $(\bl^{-n}, \bh^{-n})(s, \bX_s^n)$ up to stage $m$ is $\mc{O}(m^2)$ and the memory complexity is $\MCO(m)$.
This is infeasible in practice, as for real problems it hundreds of stages are needed. To overcome this significant problem, we propose an enhanced version of the original algorithm which reduces the time complexity to $\mc{O}(m)$ and the memory complexity to $\mc{O}(1)$. We present this new, enhanced algorithm (with pseudocode). 

\subsection{Algorithm}\label{section:edfp}
 In order to reduce the computational complexity of evaluating $(\bl^{-n}, \bh^{-n})(s, \bX_s^n)$ in the situation when $\bl^{-n}$ or $\bh^{-n}$ explicitly appears in the minimizer in \eqref{def:alphaast}, we propose the \textit{Enhanced Deep Fictitious Play}  which parametrizes both $V^n(t, \bx)$ and policy $(\ell^n, h^n)(t,\bx)$ by NNs. For simplicity, we state the algorithm based on a generic stochastic differential game, where (possibly high-dimensional) controls are denoted by $\alpha^n(t, \bx)$ for player $n$.

In each stage of the \textit{Enhanced Deep Fictitious Play}, for each planner $n$, the loss that our algorithm aims to minimize consists of two parts: (1) the loss related to solving \eqref{eq:BSDE_forward}--\eqref{eq:BSDE_backward} and (2) the error of approximating the optimal strategy $\alpha^n$ within some hypothesis spaces. The resulted approximation $\tilde \alpha^n$ will be used in the next stage of fictitious play: 
 \begin{equation}
 \begin{split}\label{def:varitional_problem}
&\inf _{Y_{0}^{n},\tilde \alpha^n,\left\{Z_{t}^{n}\right\}_{0 \leq t \leq T}}\mathbb{E}(\left|Y_{T}^{n}\right|^{2}+\tau\int_0^T \ltwonorm{\alpha^n(s,\bX_s^n)-\tilde{\alpha}^n(s,\bX_s^n)}^2 \ud s) \\
\text { s.t. } &\bX_t^n = \bx_0 + \int_{0}^{t} \mu^n(s, \bX_s^n; \tilde\balpha^{-n}(s, \bX_s^n)) \ud s + \int_{0}^{t}\Sigma(\bX_s^n)\, \mathrm{d}\bW_s,  \\
    &Y_t^n =  Y_0^n-\int_{0}^{t}g^n(s, \bX_s^n, Z_s^n; \tilde\balpha^{-n}(s, \bX_s^n))\ud s + \int_{0}^{t}(Z_s^n)\transpose\, \mathrm{d}\bW_s,\\
    &{\alpha}^n(s,\bX_s^n)=\argmin_{\beta^n} H^n(s, \bX_s^n, (\beta^n, \tilde\balpha^{-n})(s, \bX_s^n), Z_s^n),\footnotemark
\end{split}
\end{equation}
\footnotetext{Here we have assumed that the Hamiltonian $H^n$ depends on $\nabla_{\bx} V$ through $\Sigma\transpose \nabla_{\bx} V$.}where $\|\cdot\|_2$ denotes the 2-norm, $\tilde\balpha^{-n}$ denotes the collection of approximated optimal controls from the previous stage except player $n$, and $\tau$ is a hyperparameter denoting the weight between two terms in the loss function.
As detailed in Section~\ref{section:implementation}, the hypothesis space for which we search $\tilde\alpha^n$ is characterized by another NN, in addition to the one to approximate $Y_0$ and $\left\{Z_{t}^{n}\right\}_{0 \leq t \leq T}$. 
Although representing $\tilde{\alpha}^n$ with a neural network introduces approximation errors, it allows us to efficiently access the proxy of the optimal strategy $\balpha^{-n}$ in the last stage by calling corresponding networks, instead of storing and calling all the previous strategies $\balpha^{-n, m-1}, \ldots, \balpha^{-n, 1}$ due to the recursive dependence.

Numerically we solve a discretized version of \eqref{def:varitional_problem}. 
Given a partition $\pi$ of size $N_T$ on the time interval $[0,T]$, $0<=t_0<t_1<....<t_{N_T}=T$, the algorithm reads (to ease the notation, we replace the subscript $t_k$ by $k$):
\begin{align}
&\inf _{\tiny \psi_{0} \in \mathcal{N}_{0}^{n^{\prime}},\left\{\phi_{k} \in \mathcal{N}_{k}^{n}, \xi_k \in \mathcal{N}_{k}^{n^{\prime\prime}}\right\}_{k=0}^{N_{T}-1}} 
\mathbb{E}\{|Y_T^{n, \pi}|^{2}+\tau \sum_k \ltwonorm{\alpha_k^{n,\pi}-\Tilde{\alpha}_k^{n,\pi}(\bX_k^{n,\pi})}^2\Delta t_k\} \label{eq:cost_discrete}\\[0.5em] 
\text {s.t. }& \bX_{0}^{n, \pi}=\bX_{0}, \quad Y_{0}^{n, \pi}=\psi_{0}\left(\bX_{0}^{n, \pi}\right), \quad Z_{k}^{n, \pi}=\phi_{k}\left(\bX_{k}^{n, \pi}\right), \quad \tilde\alpha^{n,\pi}_k(\bX_k^{n,\pi}) = \xi_k(\bX_k^{n,\pi}),\\
&{\alpha}_k^{n,\pi} =\argmin_{\beta^n} H^{n}(t_k, \bX_{k}^{n, \pi}, (\beta^n, \tilde\balpha^{-n, \pi}_k)(\bX_{k}^{n, \pi}), Z_{k}^{n, \pi}),\quad k=0, \ldots, N_{T}-1 \\
&\bX_{k+1}^{n, \pi}=\bX_{k}^{n, \pi}+\mu^{n}\left(t_{k}, \bX_{k}^{n, \pi} ; \tilde\balpha^{-n, \pi}_k(\bX_{k}^{n, \pi})\right) \Delta t_{k}+\Sigma\left(t_{k}, \bX_{k}^{n, \pi}\right) \Delta \bW_{k},\label{eq:discrete_x}\\
&Y_{k+1}^{n, \pi}=Y_{k}^{n, \pi}-g^{n}\left(t_{k}, \bX_{k}^{n, \pi}, Z_{k}^{n, \pi} ; \tilde\balpha^{n,\pi}_k( \bX_{k}^{n, \pi})\right) \Delta t_{k}+\left(Z_{k}^{n, \pi}\right)^{\mathrm{T}} \Delta \bW_{k},\label{eq:discrete_y}
\end{align}
where $\Delta t_k=t_{k+1}-t_k$, $\Delta \bW_k=\bW_{t_{k+1}}-\bW_{t_k}$, and  $\mathcal{N}_{0}^{n}$, $\{\mathcal{N}_{k}^{n^{\prime}}\}_{k=0}^{N_T-1}$, $\{\mathcal{N}_{k}^{n^{\prime\prime}}\}_{k=0}^{N_T-1}$ are hypothesis spaces for player $n$, which will be specified later through neural network structures.

The expectation in \eqref{eq:cost_discrete} is further approximated by Monte Carlo samples of \eqref{eq:discrete_x}-\eqref{eq:discrete_y}. The parameters in the hypothesis spaces are determined by stochastic gradient descent (SGD) algorithms such that the approximated expectation is minimized, which in turn gives the optimal deterministic functions $(\psi_0^\ast, \phi_k^\ast, \xi_k^\ast)$. We expect that $(\psi_0^\ast, \phi_k^\ast, \xi_k^\ast)$ will approximate $(V^n, \nabla_\bx V^n, \alpha^n)$ well when this proxy of \eqref{eq:cost_discrete} is small. Particularly, $\{\xi^\ast_{k}\}_{k=0}^{N_T-1}$ serves as an efficient tool to evaluate the optimal policy at the current stage for finding Nash equilibrium. Implementation details and the full algorithm are presented in Section \ref{section:implementation}.
Note that when $\tau=0$ and $\tilde \alpha$ are replaced by $\alpha$ in the above algorithm, the \textit{Enhanced Deep Fictitious Play} degenerates to the \textit{Deep Fictitious Play} proposed in \cite{HaHu:19}. 

\subsection{Implementation}\label{section:implementation}
Here we provide some detail to implement the methodology in Section \ref{section:edfp}. First, we specify the hypothesis spaces for neural networks $\mathcal{N}_{0}^{n^{\prime}}$,$\{\mathcal{N}_{k}^{n}\}_{k=0}^{N_T-1}$, $\{\mathcal{N}_{k}^{n^{\prime\prime}}\}_{k=0}^{N_T-1}$, corresponding to $V^n, \nabla_\bx V^n, \alpha^n$ (the superscript $m$ is dropped again for simplicity). 
$V^n(t,x)$ is parametrized directly by a neural network $\mathrm{NN}(t,\bx)$ . Corresponding map $\Sigma(\bX)\transpose\nabla_{\bx} V^n(t, \bX)$ that defines $Z_t^n$ in the optimization problem \eqref{def:varitional_problem} could be parametrized by 
$\Sigma(\bx)\nabla_{\bx}\mathrm{NN}(t,\bx)$. Naturally, $\Sigma(\bx)\nabla_{\bx}\mathrm{NN}(t_k,\bx)$ is a hypothesis function in $\mathcal{N}_{k}^{n^{\prime}}$. Under this parametrization rule, the hypothesis functions in $\mathcal{N}_{0}^{n^{\prime}}$ and $\{\mathcal{N}_{k}^{n}\}_{k=0}^{N_T-1}$ share the same set of parameters. The policy function $\alpha^n(t,\bx)$ is parametrized by another neural network $\widetilde{\mathrm{NN}}(t,\bx)$ and then $\widetilde{\mathrm{NN}}(t_k,\bx)$ plays the role of a hypothesis function in $\mathcal{N}_{k}^{n^{\prime\prime}}$. In other words, $\{\mathcal{N}_{k}^{n^{\prime\prime}}\}_{k=0}^{N_T-1}$ share the same set of neural networks. In a stochastic game of $N$ players, there are $2N$ neural networks in total, with $N$ neural networks corresponding to $V^n(t,\bx)$ and $N$ neural networks corresponding to $\alpha ^n(t,\bx)$. At stage $m$, the $N$ $V$-networks are trained to approximate the solution of PDE \eqref{def:PDE-FK} and the $N$ $\alpha$-networks are trained to approximate the current optimal policy computed by \eqref{def:alphaast} using the optimal strategies in the last stage. The updated neural networks at stage $m$ would be used at stage $m+1$ to simulate paths $\{X_k^{n,\pi}\}_{k=0}^{N_T-1}$ and optimal strategies by \eqref{def:alphaast}. In this work, fully connected neural networks with three hidden layers are used.

Second, at each stage, the $2N$ neural networks could be decoupled to $N$ pairs of $V$-network and $\alpha$-network based on players. Then, the $N$ pairs of neural networks could be trained in parallel, which dramatically reduces computational time. As \cite{HaHu:19} and \cite{seale2006solving} pointed out, it is not necessary to solve the individual control problem accurately in each stage; the parameters at each stage are updated starting from the optimal parameters in the last stage without re-initialization. This requires only a moderate number of epochs for the stochastic gradient descent at each stage. 

The full implementation of \textit{Enhanced Deep Fictitious Play} is shown in Algorithm 1. For simplicity, we state the algorithm based on a generic stochastic differential game. 
\begin{algorithm}[!ht]
\caption{Enhanced Deep Fictitious Play for Finding Markovian Nash Equilibrium \label{def_algorithm1}}
    \begin{algorithmic}[1]
	\REQUIRE $N$ = \# of players, $N_T$ = \# of subintervals on $[0,T]$, $M$ = \# of total stages in fictitious play, $N_{\text{sample}}$ = \# of sample paths generated for each player at each stage of fictitious play, $N_{\text{SGD\_per\_stage}}$ = \# of SGD steps for each player at each stage, $N_{\text{batch}}$ = batch size per SGD update, $\balpha^0\colon$ the initial policies that are smooth enough
	    \STATE Initialize $N$ deep neural networks to represent $V^{n,0}$ and $N$ deep neural networks to represent $\alpha^{n,0}, n \in \mc{N}$
		\FOR{$m \gets 1$ to $M$}
		\FORALLP{$n \in \mc{N}$}
		\STATE  Generate $N_\text{sample}$ sample paths $\{\bX_{k}^{n,\pi}\}_{k=0}^{N_T}$ according to \eqref{eq:discrete_x} and the realized approximate optimal policies $\tilde\balpha^{-n, m-1}(t_k, \bm X_k^{n,\pi})$ 
		(Remark: $\tilde\balpha$ represents social and health policies in the multi-region SEIR model)
		\FOR{$e \gets 1$ to $N_{\text{SGD}\_\text{per}\_\text{stage}}$}
		    \STATE Update the parameters of the $n^{th}$ $V$-neural network and $\alpha$-neural network one step with $N_{\text{batch}}$ paths using the SGD algorithm (or its variant), based on the loss function \eqref{eq:cost_discrete}
		  \ENDFOR
		  \STATE Obtain the approximate optimal policy  $\tilde\alpha^{n,m}$ represented by the latest policy neural network
		  \ENDFOR
		  \STATE Collect the approximate optimal policies at stage $m$: $\tilde{\bm{\alpha}}^m \gets (\tilde\alpha^{1,m}, \ldots, \tilde\alpha^{N,m})$
		  \ENDFOR
		\RETURN The approximate optimal policy $\tilde\balpha^{M}$
	\end{algorithmic}
\end{algorithm}

Due to page limits, the exact choice of NN architectures will be detailed in Appendix~\ref{app:implementation}. To determine the total stages of fictitious play $M$, we monitor the relative changes of $\alpha^n$ and $V^n$, and stop the process when the relative change from stage to stage is below a threshold. Regarding the total number of SGD per stage,  as shown in \cite[Figure 1]{HaHu:19}, the original DFP is insensitive to the choice of $N_{\text{SGD\_per\_stage}}$. We find the enhanced version sharing the same behavior when apply to the COVID-19 case study. We give more details in Section~\ref{sec:casestudy}, and further experiments regarding different choices of $M$ and $N_{\text{SGD\_per\_stage}}$ in Appendix~\ref{app:morenumerics}.

For problems without analytical solutions, one natural concern 
	is the reliability of numerical solutions. Theoretically, the quantity \eqref{eq:cost_discrete} serves as the indicator of the numerical accuracy. In the original DFP where the second term in \eqref{eq:cost_discrete} does not exist, Theorem 3 in \cite{hanhulong:21} ensures the convergence to the true Nash equilibrium under technical assumptions when \eqref{eq:cost_discrete} is small enough for each fictitious play stage and with sufficiently large $M$ and small $\Delta t_k$. In practice, the quantity in \eqref{eq:cost_discrete} is approximated by its Monte Carlo counterpart, which we define as the loss function of our algorithms. Therefore, having small training losses during all stages will ensure convergence.
	Extending Theorem 3 in \cite{hanhulong:21} to the current setting is beyond the scope of this paper and is left for further work. 

\section{Application on COVID-19}\label{sec:case}

Our case study is based on COVID-19. We focus mainly on the lockdown/travel ban policy between different regions. Therefore, to simplify the presentation, we omit the health policy $h$ in the following discussion and make $v(\cdot) = v$, $\lambda(\cdot) = \lambda$, and $\eta = 0$. Moreover, as vaccines are not available to the population yet, we let $v(\cdot) = v = 0$.

\subsection{Parameter choices}\label{sec:parameter}

In single-region SEIR models, the transmission rate, $\beta$, is the basic reproductive number divided by the length of time an individual is infectious. In our model, we assume that there is a region-independent constant $\beta$ that underlies the rate of infections for each population. The transmission rates $\beta^{nk}$ between regions are related to the underlying transmission rate $\beta$, and the amount of travel between regions $n$ and $k$.

To quantify the size of travel between regions, we assume there is a constant fraction of people from region $n$ that travel to region $k$, $f^{nk}$, at any given moment in time. We note that realistically one may expect $f^{nk}$ to depend on time and also on the epidemic status of regions $n$ and $k$. However, for simplicity, we will not consider these scenarios in our numerical experiments. We assume that $f^{nn} \gg f^{nk}$, $f^{nn} \gg f^{kn}$, $ \forall k \neq n$, meaning that most of the population $n$ resides in region $n$ at any given time, and also that most of the people in region $n$ at a given time are from $n$ and not travelers from another region. We will see later that this implies $\beta^{nn} \gg \beta^{nk}$ $ \forall k \neq n$.

To further clarify the transmission of infection from one region due to another, we need to describe the set of parameters $\{ \beta^{nk} : n,k \in \cN \}$. Here, $\beta^{nk}$ represents the rate of transmission from region $k$ to region $n$. Specifically, $\beta^{nk}$ is the number of infected people in population $n$ per a contactable, infectious individual in population $k$ per day,  assuming that $100\%$ of population $n$ is susceptible. This definition of $\beta^{nk}$ comes from the derivation of the SDE system itself \eqref{def:St}--\eqref{def:It}. Accounting for infection in region $n$ by individuals in region $k$ from travel to both region $n$ and $k$, we  have that

\begin{equation}\label{def:Beta}
\beta^{nk} = \begin{dcases}
 \beta (f^{nk} f^{kk} + f^{kn}f^{nn}) \frac{P^k}{P^n},  & \mathrm{if} \  k \neq n \\
 \beta (f^{nn})^2,  & \mathrm{if} \  k = n.\\
\end{dcases}
\end{equation}
We defer the detailed derivation of \eqref{def:Beta} to Appendix~\ref{app:beta}.

Therefore, to specify $\beta^{nk}$, we need to provide $\beta$ and $f^{nk}$. We will specify $f^{nk}$ for New York (NY), New Jersey (NJ), and Pennsylvania (PA) in the next section. To estimate $\beta$, we choose a basic reproductive number $R_0 = 2.2$, which is consistent with \cite{fauci}, and assume that the length of each individual being infectious is $13$ days. More precisely, we assume infected individuals either recover or die in 13 days. Under these assumptions, we obtain $\beta = \frac{2.2}{13} \approx 0.17$, consistent with \cite{deaths} as a 13 day median time until death from illness onset is used.

The infection fatality rate, or IFR, is the fraction of those infected who died from the infection. We choose the IFR to be 0.65\% according to the CDC estimate. This is also consistent with \cite{Meyerowitz-Katz2020.05.03.20089854}, which suggests a point estimate of 0.68\%. 
The assumptions of an IFR of 0.65\% and an infectious period of 13 days determines that the recovery rate (including both recovery and death due to infection) is $\lambda =  \frac{1}{13} \approx 0.0769$, and the death rate is $\kappa =  \frac{(0.65\%)}{13} = 0.0005$. We choose the latent period to be 5 days according to \cite{Incubation}. This means that the we will have $\gamma = \frac{1}{5}$. Note that this choice has also been used in other models such as \cite{econpaper} and \cite{peng2020epidemic}. We assume that the parameters for noise-level $\sigma_{s_n}, \sigma_{e_n}$, $n \in \cN$ are all $0.0002$, and the extent to which one adheres to the social distancing policy, $\theta$, is either $\theta = 0.9$ or $\theta = 0.99$. 

With most of the parameters for the SDE model \eqref{def:St}-\eqref{def:It} discussed, we now address those specific to defining the cost. Regarding the risk-free-rate $r$, note that U.S. Treasury yields are historically low and the uncertainty in the current level of inflation. We choose for simplicity that $r=0$. Also, considering that we are interested in simulations with time periods of less than a year, the discounting is negligible. The parameter $w$ represents the dollar output per individual per day. To estimate $w$, we use GDP per capita per day, yielding the estimate $w = 172.6$ dollars per person per day. Following \cite{econpaper} and \cite{othereconpaper}, we use the value of a statistical life, $\chi$, to be 20 times GDP per capita. This results in  $\chi = 1.95 \cdot 10^6$  dollars per person. According to the CDC summary of U.S. COVID-19 activity, the hospitalization rate was 228.7 per 100,000 population by 11/14/2020. Thus, we set $p = 228.7\times 10^{-5}$. The cost per inpatient day is $c = 73300/13$ dollars, estimated according to \cite{inpatientcost}. The attention hyperparameter $a$ takes various values in the case study, and will be specified in Section~\ref{sec:casestudy}.

\subsection{NY-NJ-PA COVID-19 case study}\label{sec:casestudy}

In this section, we apply our model \eqref{def:St}--\eqref{def:J} to analyze COVID-19 related policy in three adjacent states:  New York, New Jersey, and Pennsylvania. This case study is done over 180 days starting from 03/15/2020, using the Enhanced Deep Fictitious Play algorithm introduced in Section~\ref{sec:algorithm} (cf. Algorithm~\ref{def_algorithm1}) and the parameters discussed in~\ref{sec:parameter}. The exact formulas of $\mu^n$ and $h^n$ in equation \eqref{def:PDE-FK} are derived in Appendix~\ref{app:details-sec3}.

We refer to New York State as region 1, New Jersey State as region 2, and Pennsylvania State as region 3. Their respective populations are $P^1 = 19.54$ million, $P^2 = 8.91$ million, and $P^3 = 12.81$ million. Regarding $\beta^{nk}$, $\forall n, k = 1, 2, 3$, we assume that: (a) 90\% of any state's population is residing in their state at a given time; (b) the remaining population (travelers) visit the other regions in an equal proportion; and (c) there is no travel outside of the considered regions, {\it i.e.}, the NY-NJ-PA is a closed system. The reasoning for (c) is that, under our model assuming that infection only occurs in the regions considered, (c) is equivalent to allowing people traveling outside the considered regions, but the travelers cannot be affected. For simplicity, we assume this is the case. Under these assumptions, we will have $f^{nn} = 90\%$ for $n=1,2,3$ and $f^{nk} = 5\%$ for $n \neq k$, and obtain the values of $\beta^{nk}$ through \eqref{def:Beta}.

\begin{figure}[h]
    \centering
    \includegraphics[width = 0.75\textwidth, trim = {2em 2em 5em 6em}, clip, keepaspectratio=True]{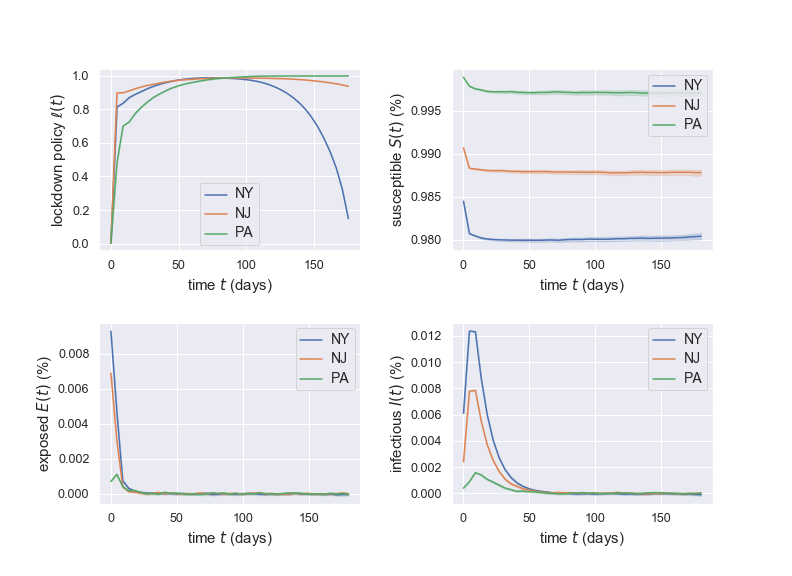}
    \caption{Plots of optimal policies (top-left), Susceptibles (top-right), Exposed (bottom-left) and Infectious (bottom-right) for three states: New York (blue), New Jersey (orange) and Pennsylvania (green). The shaded areas depict the mean and 95\% confidence interval over 256 sample paths. Choices of parameters are in Section~\ref{sec:parameter}, $a = 100$ and $\theta = 0.99$.}
    \label{fig:control}
\end{figure}

\begin{figure}[h]
    \centering
    \includegraphics[width = 0.75\textwidth, trim = {2em 2em 5em 6em}, clip, keepaspectratio=True]{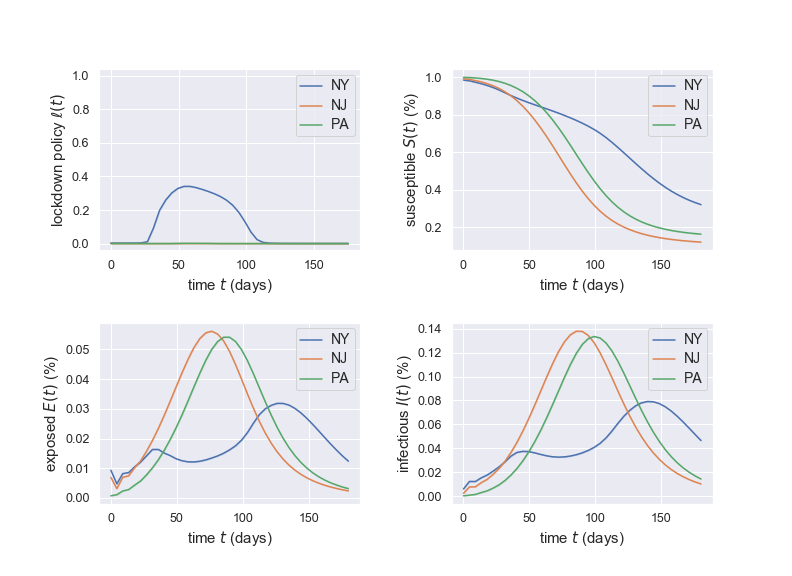}
    \caption{An illustration that governors' inaction or mild control leads to disease spreading.}
    \label{fig:outofcontrol}
\end{figure}

Figure~\ref{fig:control} presents the equilibrium policy issued by the governors of NY, NJ, and PA  when the policy effectiveness is $\theta = 0.99$, {\it i.e.}, 99\% of the population follow the lockdown order. The hyperparameter is $a = 100$, {\it i.e.}, each governor values people's death 100 times the lockdown cost. In this scenario, the governors take action at an early stage and soon reach the strictest policy. Once the disease is under control, they may relax the policy later. The percentage of Susceptible, Exposed, Infectious, and Removed stays almost constant in the end. As a comparison, Figure~\ref{fig:outofcontrol} illustrates how the pandemic gets out of control if governors show inaction or issue mild lockdown policies.

\smallskip

\noindent\textbf{Experiment 1: dependence on $a$.}
We further analyze how the planners' view on 
the death of human beings changes their policies. In reality, economic loss is not the only factor the planners concern about. It is also important to mitigate the infections and deaths within the budget and available resources. Different views and values from the planners will lead to different policies. In this experiment, we consider different attitudes towards the infection, especially death caused by COVID-19. This is reflected by the attention hyperparameter $a$. Large $a$ implies that planners care more about human beings and are willing to spend more effort or endure more economic loss on lockdown to avoid further infection and death. In comparison, smaller $a$ implies that planners care less about  infection and death and instead pay more attention to minimizing the total cost. 

\begin{figure}[ht]
    \centering
    \includegraphics[width = 0.75\textwidth, trim = {2em 2em 5em 4em}, clip, keepaspectratio = True]{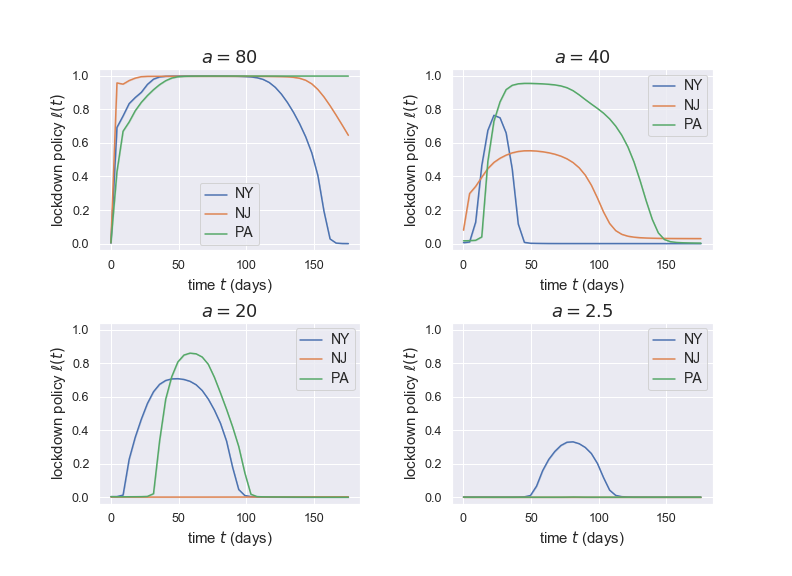}
    \caption{Plots of optimal policies with different choice of $a$ for three states: New York (blue), New Jersey (orange) and Pennsylvania (green), when the lockdown efficiency is $\theta = 0.9$.}
    \label{fig:theta09}
\end{figure}

\begin{figure}[ht]
    \centering
    \includegraphics[width = 0.75\textwidth, trim = {2em 2em 5em 4em}, clip, keepaspectratio = True]{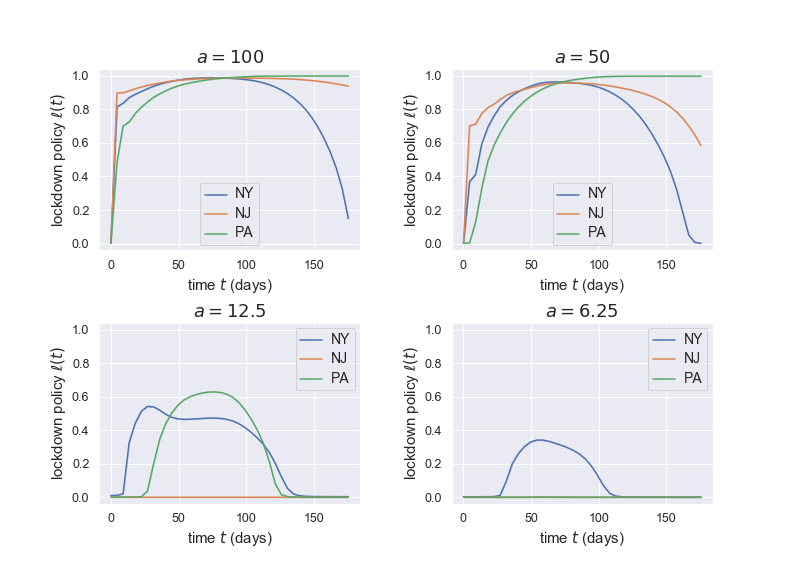}
    \caption{Plots of optimal policies with different choice of $a$ for three states: New York (blue), New Jersey (orange) and Pennsylvania (green), when the lockdown efficiency is $\theta = 0.99$.}
    \label{fig:theta099}
\end{figure}

The numerical results in Figures \ref{fig:theta09} and \ref{fig:theta099} are consistent with intuition. With a large $a$ (top-left panels), meaning the planners give more consideration to infection and death, they tend to issue a restrict lockdown policy, which helps slow down the disease spreads and reduce the percentage of infected people. As $a$ becomes smaller (top-right panels), planners weigh more the economic loss and spend fewer efforts on lockdown.  
When the attention $a$ is small enough, some states even give up controlling the disease spread due to economic concern (bottom panels). As a result, the pandemic would get out of control by the end of the simulation period. This mild lockdown policy leads to a natural spread of disease (also shown in Figure \ref{fig:outofcontrol}).

\smallskip

\noindent\textbf{Experiment 2: dependence on $\theta$.} We next analyze how the residents' willingness to comply with the lockdown policy changes the optimal policies and the development of a pandemic. The larger the $\theta$ is, the more likely the residents will follow the lockdown policy, and the larger the difference the control makes on the pandemic situation. Conversely, small $\theta$ weakens the effect of the lockdown policy. In the extreme case of $\theta=0$, no matter how strict the lockdown policy is, the pandemic will become a natural spread because the control term in \eqref{def:St} disappears. In short, this willingness to policy compliance should be an essential factor in decision-making. 

\begin{figure}[!ht]
    \centering
    \includegraphics[width = 0.75\textwidth, trim = {3em 6em 6em 10em}, clip, keepaspectratio = True]{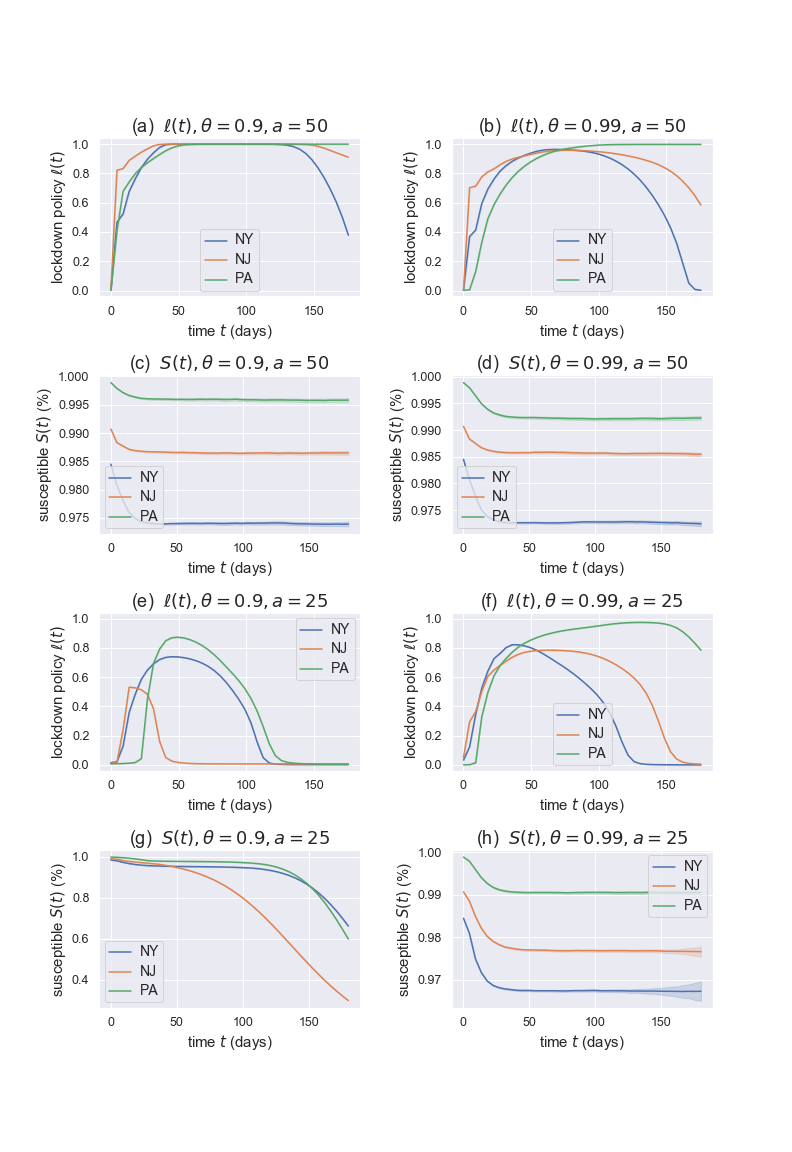}
    \vspace{-2.4em}
    \caption{Comparison of optimal policies for three states (NY = blue, NJ = orange, PA = green) and their susceptibles between different policy effectivenss $\theta$ and hyperparameter $a$.}
    \label{fig:comparison}
\end{figure}

To this end, we compare the optimal policy when $\theta=0.9$ and $\theta=0.99$ in Figure \ref{fig:comparison}. Panels (a-d) show the difference of optimal policies $\ell(t)$ and the Susceptible $S(t)$ in the tri-state game under different $\theta$ when $a=50$. In both situations, the pandemic is well-controlled, with the percentage of susceptible people staying stable in the end. Moreover, in the case of $\theta = 0.99$, people are more willing to comply with the policies. Consequently,  the planners are allowed to use a less strict lockdown policy as shown in Figure \ref{fig:comparison}(b) compared to \ref{fig:comparison}(a), which saves the lockdown cost. Figure \ref{fig:comparison} (e-h) shows an interesting case in the comparison of $\theta=0.9$ and $\theta=0.99$. In this scenario, with the same attention parameter ($a=25$), $\theta = 0.9$ leads to a mild lockdown policy, see Figure \ref{fig:comparison}(e), while $\theta = 0.99$ provides a possibility to stop the spread of virus, see Figure \ref{fig:comparison}(f). We believe that the decision when $\theta = 0.9$ is a compromise as the lockdown is not efficient enough to reduce largely the infection and death loss by paying lockdown cost, and also due to the limited simulation period, {\it i.e.}, the policies could have been different if we had the simulation until the disease dies out. We also believe that the early give-up by NJ drives NY and PA to lift lockdown policies at a later stage, because even NY and PA issue strict policies, they are still facing severe infections from NJ due to its high infected percentages and the existence of travel between states whatever the policy is. So their interventions are not worth the candle. 
Figure \ref{fig:comparison}(f)(h) further elucidate the importance of residents' support in slowing down the pandemic. Further experiments based on different sets of $(a, \theta)$ reveal the possibility of having multiple Nash equilibrium, with more elaboration in Appendix~\ref{app:NE}.

To summarize, the numerical experiments illustrate that both the balance of economy and infection/death from the view of plan-makers and the willingness of residents to follow the lockdown policy play an important role in decision-making. In reality, all three states issued stay-at-home orders in March, and attempted to reopen in June. By comparing real world policies and our simulations of $\bl(x)$, we may infer $\alpha$ and $\theta$ for NY, NJ, and PA in our model, {\it i.e.}, $\theta = 0.99$ and $a = 25$.

\section{Conclusion}\label{sec:conclusion}

In this paper, we propose a novel multi-region SEIR model to study optimal policies under a  pandemic. Our new model, built on game theory, takes into account how the social and health policies issued by multiple region planners affect the progress of infectious diseases. This feature makes the model more realistic and powerful but also introduces a formidable computational challenge due to the high-dimensionality of the solution space and the strong coupling of planners' policies. We propose the {enhanced deep fictitious play} algorithm to overcome the curse of dimensionality and use the model and algorithm in a case study of the COVID-19 pandemic in three states, New York, New Jersey, and Pennsylvania. The model parameters are estimated from  real data posted by the CDC. We are able to show the effect of lockdown/travel ban policy on the spread of COVID-19 for each state and how people's willingness to comply and planners' attitude towards deaths influence the equilibrium strategies as a consequence of the competition between regions. We hope our model can draw more attention to studying optimal interventions in infectious diseases using game theory. Our numerical simulations can shed light on public policies. 

In reality, during a pandemic, the planning is usually for short periods and adjusted frequently. This can be modeled using repeated games, and planners may infer other regions' cost functional from past game outcomes. The assumptions that some parameters are identical across different regions can be relaxed and the health policy can also be added for more accurate simulations. These will be left for future work.

\acks{R.H.~was partially supported by  NSF grant DMS-1953035. H.D.C.~acknowledges partial support from NSF grant DMS-1818821. This project was jointly supervised by R.H. and H.D.C.. Y.X. and R.B. have equal contribution as the first authors. J.H. contributes on the algorithms and early discussions of the model setup.}

\pagebreak
\bibliographystyle{plain}
\bibliography{references}

\begin{thebibliography}{25}
\providecommand{\natexlab}[1]{#1}
\providecommand{\url}[1]{\texttt{#1}}
\expandafter\ifx\csname urlstyle\endcsname\relax
  \providecommand{\doi}[1]{doi: #1}\else
  \providecommand{\doi}{doi: \begingroup \urlstyle{rm}\Url}\fi

\bibitem[Alvarez et~al.(2020)Alvarez, Argente, and Lippi]{econpaper}
F.~E. Alvarez, D.~Argente, and F.~Lippi.
\newblock A simple planning problem for {COVID}-19 lockdown.
\newblock Working Paper 26981, National Bureau of Economic Research, 2020.

\bibitem[Bauch and Earn(2004)]{bauch2004vaccination}
C.~T. Bauch and D.~J.~D. Earn.
\newblock Vaccination and the theory of games.
\newblock \emph{Proceedings of the National Academy of Sciences}, 101\penalty0
  (36):\penalty0 13391--13394, 2004.

\bibitem[Bauch et~al.(2003)Bauch, Galvani, and Earn]{bauch2003group}
C.~T. Bauch, A.~P. Galvani, and D.~J.~D. Earn.
\newblock Group interest versus self-interest in smallpox vaccination policy.
\newblock \emph{Proceedings of the National Academy of Sciences}, 100\penalty0
  (18):\penalty0 10564--10567, 2003.

\bibitem[Brown(1949)]{Br:49}
G.~W. Brown.
\newblock Some notes on computation of games solutions.
\newblock Technical report, Rand Corp Santa Monica CA, 1949.

\bibitem[Brown(1951)]{Br:51}
G.~W. Brown.
\newblock Iterative solution of games by fictitious play.
\newblock \emph{Activity Analysis of Production and Allocation}, 13\penalty0
  (1):\penalty0 374--376, 1951.

\bibitem[Chang et~al.(2020)Chang, Piraveenan, Pattison, and
  Prokopenko]{gamereviewpandemic}
S.~L. Chang, M.~Piraveenan, P.~Pattison, and M.~Prokopenko.
\newblock {Game theoretic modelling of infectious disease dynamics and
  intervention methods: a review}.
\newblock \emph{Journal of Biological Dynamics}, 14\penalty0 (1):\penalty0
  1--33, 2020.

\bibitem[E et~al.(2017)E, Han, and Jentzen]{EHaJe:17}
W.~E, J.~Han, and A.~Jentzen.
\newblock Deep learning-based numerical methods for high-dimensional parabolic
  partial differential equations and backward stochastic differential
  equations.
\newblock \emph{Communications in Mathematics and Statistics}, 5\penalty0
  (4):\penalty0 349--380, 2017.

\bibitem[El~Karoui et~al.(1997)El~Karoui, Peng, and Quenez]{ElPeQu:97}
N.~El~Karoui, S.~Peng, and M.~C. Quenez.
\newblock Backward stochastic differential equations in finance.
\newblock \emph{Mathematical Finance}, 7\penalty0 (1):\penalty0 1--71, 1997.

\bibitem[Fauci et~al.(2020)Fauci, Lane, and Redfield]{fauci}
A.~S. Fauci, H.~C. Lane, and R.~R. Redfield.
\newblock {COVID}-19 —navigating the uncharted.
\newblock \emph{New England Journal of Medicine}, 382\penalty0 (13):\penalty0
  1268--1269, 2020.

\bibitem[Hall et~al.(2020)Hall, Jones, and Klenow]{othereconpaper}
R.~E. Hall, C.~I Jones, and P.~J. Klenow.
\newblock Trading off consumption and {COVID}-19 deaths.
\newblock Working Paper 27340, National Bureau of Economic Research, 2020.

\bibitem[Han and Hu(2020)]{HaHu:19}
J.~Han and R.~Hu.
\newblock Deep fictitious play for finding {Markovian} {Nash} equilibrium in
  multi-agent games.
\newblock In \emph{Proceedings of The First Mathematical and Scientific Machine
  Learning Conference (MSML)}, volume 107, pages 221--245, 2020.

\bibitem[Han et~al.(2018)Han, Jentzen, and E]{HaJeE:18}
J.~Han, A.~Jentzen, and W.~E.
\newblock Solving high-dimensional partial differential equations using deep
  learning.
\newblock \emph{Proceedings of the National Academy of Sciences}, 115\penalty0
  (34):\penalty0 8505--8510, 2018.

\bibitem[Han et~al.(2020{\natexlab{a}})Han, Hu, and Long]{han2020convergence}
J.~Han, R.~Hu, and J.~Long.
\newblock Convergence of deep fictitious play for stochastic differential
  games.
\newblock \emph{arXiv preprint arXiv:2008.05519}, 2020{\natexlab{a}}.

\bibitem[Han et~al.(2020{\natexlab{b}})Han, Hu, and Long]{hanhulong:21}
J.~Han, R.~Hu, and J.~Long.
\newblock Barron metric for the convergence of empirical distribution.
\newblock \emph{in preparation}, 2020{\natexlab{b}}.

\bibitem[Health(2020)]{inpatientcost}
Fair Health.
\newblock {Costs for a Hospital Stay for COVID-19}, 2020.
\newblock URL
  \url{https://www.fairhealth.org/article/costs-for-a-hospital-stay-for-covid-19}.

\bibitem[Hu(2020)]{Hu2:19}
R.~Hu.
\newblock Deep fictitious play for stochastic differential games.
\newblock \emph{Communications in Mathematical Sciences}, 2020.

\bibitem[Isaacs(1965)]{Isaacs1965}
R.~Isaacs.
\newblock \emph{Differential Games: A Mathematical Theory with Applications to
  Warfare and Pursuit, Control and Optimization}.
\newblock London: John Wiley and Sons, 1965.

\bibitem[Lauer et~al.(2020)Lauer, Grantz, Bi, Jones, Zheng, Meredith, Azman,
  Reich, and Lessler]{Incubation}
S.~A. Lauer, K.~H. Grantz, Q.~Bi, F.~K. Jones, Q.~Zheng, H.~R. Meredith, A.~S.
  Azman, N.~G. Reich, and J.~Lessler.
\newblock The incubation period of coronavirus disease 2019 ({COVID}-19) from
  publicly reported confirmed cases: Estimation and application.
\newblock \emph{Annals of Internal Medicine}, 172\penalty0 (9):\penalty0
  577--582, 2020.

\bibitem[Linton et~al.(2020)Linton, Kobayashi, Yang, Hayashi, Akhmetzhanov,
  Jung, Yuan, Kinoshita, and Nishiura]{deaths}
N.M. Linton, T.~Kobayashi, Y.~Yang, K.~Hayashi, A.R. Akhmetzhanov, S.-M. Jung,
  B.~Yuan, R.~Kinoshita, and H.~Nishiura.
\newblock Incubation period and other epidemiological characteristics of 2019
  novel coronavirus infections with right truncation: A statistical analysis of
  publicly available case data.
\newblock \emph{Journal of Clinical Medicine}, 538\penalty0 (9), 2020.

\bibitem[Liu et~al.(1987)Liu, Hethcote, and Levin]{dynamicalliu}
W.-M. Liu, H.~W. Hethcote, and S.~A. Levin.
\newblock {Dynamical behavior of epidemiological models with nonlinear
  incidence rates}.
\newblock \emph{Journal of Mathematical Biology}, 25\penalty0 (4):\penalty0
  359--380, 1987.

\bibitem[Meyerowitz-Katz and Merone(2020)]{Meyerowitz-Katz2020.05.03.20089854}
G.~Meyerowitz-Katz and L.~Merone.
\newblock A systematic review and meta-analysis of published research data on
  {COVID}-19 infection-fatality rates.
\newblock \emph{medRxiv}, 2020.
\newblock \doi{10.1101/2020.05.03.20089854}.

\bibitem[Pardoux and Peng(1992)]{PaPe:92}
E.~Pardoux and S.~Peng.
\newblock Backward stochastic differential equations and quasilinear parabolic
  partial differential equations.
\newblock In \emph{Stochastic Partial Differential Equations and Their
  Applications}, pages 200--217. Springer, 1992.

\bibitem[Pardoux and Tang(1999)]{PaTa:99}
E.~Pardoux and S.~Tang.
\newblock Forward-backward stochastic differential equations and quasilinear
  parabolic {PDEs}.
\newblock \emph{Probability Theory and Related Fields}, 114\penalty0
  (2):\penalty0 123--150, 1999.

\bibitem[Peng et~al.(2020)Peng, Yang, Zhang, Zhuge, and Hong]{peng2020epidemic}
L.~Peng, W.~Yang, D.~Zhang, C.~Zhuge, and L.~Hong.
\newblock Epidemic analysis of covid-19 in china by dynamical modeling, 2020.

\bibitem[Seale and Burnett(2006)]{seale2006solving}
D.~Seale and J.~Burnett.
\newblock Solving large games with simulated fictitious play.
\newblock \emph{International Game Theory Review}, 8\penalty0 (03):\penalty0
  437--467, 2006.

\end{thebibliography}

\appendix

\section{Technical Details to the Stochastic Multi-region SEIR Model}\label{app:details}

\subsection{The dynamics of \texorpdfstring{$\bX_t$}{} in Section~\ref{sec:model}}\label{app:details-sec2}
In Section~\ref{sec:model}, for the ease of notations and clarity of the presentation, we rewrite the dynamics of $(S_t^n, E_t^n, I_t^n)$ defined in \eqref{def:St}--\eqref{def:It} in the vector form
\begin{equation}\label{app:Xt}
    \ud \bX_t = b(t, \bX_t, \bm \ell(t, \bX_t), \bm h(t, \bX_t)) \ud t  + \Sigma(\bX_t)\ud \bm W_t, 
\end{equation}
where $\bX_t \equiv [\bm S_t,\bm E_t, \bm I_t]\transpose \equiv [S_t^1, \cdots, S_t^N, E_t^1, \cdots, E_t^N, I_t^1, \cdots, I_t^N]\transpose \in \RR^{3N}$, the Markovian controls $(\bl, \bh)$ are given by $\bl(t, \bx) = [\ell^1, \ldots, \ell^N]\transpose(t,\bx)$ and  $\bh(t, \bx) = [h^1, \ldots, h^N]\transpose(t,\bx)$. In the sequel, for a vector $\bx \equiv (\bs, \be, \bi) \in \RR^{3N}$, we shall index them in two ways,
\begin{equation}
    (s_1, \cdots, s_N, e_1, \cdots, e_N, i_1, \cdots, i_N) \text{ or } (x_1, \ldots, x_{3N}).
\end{equation}
and use them interchangeably. The former one emphasizes the dependence on each category, while the later notation is more condensed. Similarly, for partial derivatives, we will have two set of notations 
\begin{equation}
    (\partial_{s_1}, \cdots, \partial_{s_N}, \partial_{e_1}, \cdots, \partial_{e_N}, \partial_{i_1}, \cdots, \partial_{i_N}) \text{ or } (\partial_{x_1}, \cdots, \partial_{x_{3N}}).
\end{equation}

We give precise definitions of \eqref{app:Xt} in this appendix.  $b(t, \bx, \bl, \bh) = [b_1, \ldots, b_{3N}]\transpose(t, \bx, \bl, \bh)$ is a deterministic vector-valued function:
\begin{equation}
b_j (t, \bx, \bl, \bh) =
\begin{dcases}
       -\sum_{k = 1}^N \beta^{jk} s_ji_k (1-\theta \ell^j(t, \bx))  (1-\theta \ell^k(t, \bx)) - v(h^j(t, \bx))s_j,  \qquad   j \in \cN, \\
    \sum_{k = 1}^N \beta^{\tj k} s_\tj i_k (1-\theta \ell^\tj(t, \bx))  (1-\theta \ell^k(t, \bx)) - \gamma e_\tj,  \qquad \quad j \in \cN + N,  \\
   \gamma e_\tj - \lambda(h^\tj(t, \bx)) i_\tj,  \qquad \qquad  j \in \cN + 2N, \text{ and } \tj = j~\text{mod}~N.
\end{dcases}
\end{equation}
$\Sigma(\bx) = (\Sigma_{j,k}(\bx))$ is a matrix-valued deterministic function in  $\RR^{3N \times 2N}$ with non-zero entries given below:
\begin{align}
&\Sigma_{j,j}(\bx) = -\sigma_{s_j} s_j,  \quad &&\Sigma_{j+N,j}(\bx) = \sigma_{s_j} s_j, \\
&\Sigma_{j+N,j+N}(\bx) = -\sigma_{e_j} e_j, \quad &&\Sigma_{j+2N,j+N}(\bx) = \sigma_{e_j} e_j, \quad j \in \cN.
\end{align}
and $\{\bm W_t\}_{0\leq t \leq T}$ is a $2N$-dimensional standard Brownian motion:
$$\bm W_t = [W^{s_1}_t, \cdots, W^{s_N}_t, W^{e_1}_t, \cdots, W^{e_N}_t]\transpose.$$

Each region's running cost $f^n$ defined in \eqref{def:cost} is
\begin{equation}
    f^n(t, \bx, \bl, \bh) = e^{-rt}P^n[(s_n + e_n + i_n)\ell^n(t, \bx) w + a(\kappa i_n \chi +  p i_n c)] + e^{-rt}\eta (h^n(t, \bx))^2. 
\end{equation}

\subsection{The decoupled HJB equations in the form of
\texorpdfstring{\eqref{def:PDE-FK}}{}}\label{app:details-sec3}

Recall that we aim to solve \eqref{def:HJB-DFP} using the BSDE approach (nonlinear Feynman Kac relation). To this end, in this appendix, we will rewrite it in the form of \eqref{def:PDE-FK} and identify $\mu^n$ and $g^n$. The first step is to identify the minimizer in the Hamiltonian \eqref{def:H}. Keeping in mind that our testing case is COVID-19 where vaccines are not fully developed yet, the term $v(h_t^n) = 0$ is essentially zero in the past. Also, to focus on the lockdown/travel ban policy between different regions (as we did in the case study), we shall exclude the health policy $\bh(t, \bx)$ from planners' decision problem, {\it i.e.}, $v(\cdot) = 0$, $\lambda(\cdot) = \lambda$, and $\eta = 0$ in the following derivations, and remove the dependence of $\bh$ in all relevant functions. We remark that including the health policy $\bh(t, \bx)$ is a straightforward generalization.

Recall that the Hamiltonian in \eqref{def:HJB-DFP} reads:
\begin{align}
   & H^n(t, \bx, (\ell^n, \bl^{-n,m}), \nabla_\bx V^{n, m+1}) \\
&   \qquad  = b(t, \bx, (\ell^n, \bl^{-n,m})) \cdot \nabla_\bx V^{n, m+1} + f^n(t, \bx, \ell^n)\\
 &   \qquad  = \sum_{j=1}^{3N}b_j(t, \bx, (\ell^n, \bl^{-n,m})) \pdv{V^{n, m+1}}{x_j} + e^{-rt}P^n[(s_n + e_n + i_n)\ell^n(t, \bx) w + a(\kappa i_n \chi +  p i_n c)],
\end{align}
and recall that $\bl^{-n,m} = (\ell^{1, m}, \ldots,  \ell^{n-1, m}, \ell^{n+1, m}, \ldots, \ell^{N,m})$ represents the $m^{th}$ stage strategies of all players other than $n$, which are given functions in this derivation. The first order condition requires for $\ell^n$:
\begin{align}
      0 = \sum_{\substack{j = 1\\ j\neq n}}^N (1-\theta\ell^{j,m})\left[\dm{j}{n}{n}{m+1} + \dm{n}{j}{n}{m+1}\right]  \\
    + 2(1-\theta\ell^n) \dm{n}{n}{n}{m+1} -\frac{1}{\theta} e^{-rt}P^n(s_n + e_n + i_n)w.
\end{align}
The critical point given by the above equation indeed gives a minimizer of the Hamiltonian, as long as it is in $[0,1]$. Because we can show $\left(\pdv{V^{n, m+1}}{e_n} - \pdv{V^{n, m+1}}{s_n}\right) > 0$ by comparing $V^{n, m+1}(t, \bx + \eps_{n+N})$ and $V^{n, m+1}(t, \bx + \eps_n)$ using their definitions, where $\eps_{j}$ is a $3N$-vector with only one nonzero entry $\eps \ll 1$ at $j^{th}$ position. Intuitively, with all others players' initial condition the same, if player $n$ starts with a higher exposed proportion $e_n + \eps$, it will produce more cost, comparing with the same increase proportion still being susceptible $s_n + \eps$. To summarize, we deduce the optimal policy for player $n$ at stage $m+1$ is given by:
\begin{align}\label{def:alpha-DFP}
    \ell^{n, m+1}(t, \bx) &=  \Bigg\{2\dm{n}{n}{n}{m+1} - \frac{1}{\theta} e^{-rt}P^n(s_n + e_n + i_n)w \\
 +\sum_{\substack{j = 1\\ j\neq n}}(1-&\theta \ell^{j, m})\left[\dm{j}{n}{n}{m+1} + \dm{n}{j}{n}{m+1}\right] \Bigg\}\\
& \times \left\{2\theta \dm{n}{n}{n}{m+1}\right\}^{-1}\wedge 1 \vee 0,
\end{align} 
where we use the conventional notations $a \wedge b = \min\{a, b\}$ and $a \vee b = \max\{a, b\}$.
Plugging \eqref{def:alpha-DFP} into equation \eqref{def:HJB-DFP} and by straightforward computation, one obtains for the $(m+1)^{th}$ stage,  $\mu^{n, m+1} (t, \bx; \bl^{-n,m}) =  [\mu_1^{n, m+1}, \ldots, \mu_{3N}^{n, m+1}](t, \bx; \bl^{-n,m})\transpose$ is
\begin{align}
    & \mu_j^{n, m+1} = -\beta^{jn}s_ji_n (1-\theta \ell^{j,m}(t, \bx)) - \sum_{\substack{k = 1\\ k\neq n}}^N \beta^{jk}s_ji_k(1-\theta\ell^{j,m}(t, \bx))(1-\theta\ell^{k,m}(t, \bx)), \\
    &  \hspace{200pt} j \in \cN \setminus n\\
    & \mu_{n}^{n, m+1} = - \beta^{nn} s_n i_n - \sum_{\substack{k = 1\\ k\neq n}} \beta^{nk} s_n i_k (1-\theta \ell^{k, m}(t, \bx)) \\
    & \mu_{N+j}^{n, m+1} = - \mu_j^{n, m+1} - \gamma e_j, \quad \mu_{2N+j}^{n, m+1} = \gamma e_j - \lambda i_j, \quad j \in \cN.
\end{align}
To write $g^{n, m+1}$ as a function of $(t, \bx, z)$, we first compute 
{\small
\begin{align}
    &\Sigma(\bx)\transpose \nabla_\bx V^{n}(t, \bx) \\
    &=  \left[ \sigma_{s_1} s_1 \big(\pdv{V^n}{e_1} - \pdv{V^n}{s_1}\big), \cdots, \sigma_{s_N} s_N \big(\pdv{V^n}{e_N} - \pdv{V^n}{s_N}\big), \sigma_{e_1} e_1 \big(\pdv{V^n}{i_1} - \pdv{V^n}{e_1}\big), \cdots , \sigma_{e_N} e_N \big(\pdv{V^n}{i_N} - \pdv{V^n}{e_N}\big)\right]\transpose.
\end{align}
}
and then $g^{n, m+1}$ is given by:
{\small
\begin{align}
    &g^{n, m+1}(t, \bx, z; \bl^{-n,m}) = \frac{\theta^2}{\sigma_{s_n}}\beta^{nn}z_ni_n[\ell^{n, m+1}(t, \bx)]^2 \\
    & + \left\{e^{-rt}P^n(s_n + e_n + i_n)w - 2\frac{\theta}{\sigma_{s_n}}\beta^{nn}z_ni_n - \sum_{\substack{j =1 \\j \neq n}}^N \theta (1-\theta \ell^{j,m}(t, \bx))(\frac{\beta^{nj}}{\sigma_{s_n}}z_ni_j + \frac{\beta^{jn}}{\sigma_{s_j}}z_ji_n)\right\}  \ell^{n, m+1}(t, \bx)\\
    & +     e^{-rt}P^na (\kappa i_n \chi + p i_n c). \\
\end{align}
}

\subsection{The derivation of the transmission rate $\beta^{nk}$ in \texorpdfstring{\eqref{def:Beta}}{}}\label{app:beta}
We denote by $\beta$ the underlying transmission rate of the virus, which is assumed to be region independent. This transmission rate is the average number of people infected by an infectious person per day (assuming that the susceptible population is 100\%).  Thus, if a proportion $S$ of the population is susceptible, then $\beta S$ represents the average number of people infected per infectious person per day. If there are a total of $P$ people in a population with a fraction $I$ being infectious, then $\beta S (PI)$ is the number of newly infected per day. 
	
	Thus, in the context of a single-region SEIR model, the number of newly infected (or the influx to the exposed population) that occurs within $(t,t+ \ud t)$ is given by $\beta S(t) (I(t)P) \ud t$.  Dividing by $P$, the influx to the scaled exposed population in $(t,t+ \ud t)$ is $\beta S(t) I(t) \ud t$. 
	
	Now let us consider the multi-region case and temporarily ignore the effect of lockdown. The term from \eqref{def:St}--\eqref{def:It} that gives the influx to $E^n$ due to infection from $I^k$ is $\beta^{nk}S^n(t)I^k(t)\ud t$. To determine $\beta^{nk}$, we build this exact influx from core assumptions.

	First, we quantify the number of people from region $n$ that are infected by those from region $k$. Specifically, the influx in the interval $(t,t+\ud t)$ to the unscaled exposed population $n$ due to transmission from population $k$ is given by
	\begin{equation}\label{def:flux}
	\sum_{\ell} (\beta f^{n\ell} S^n(t))(f^{k\ell}I^k(t)P^k) \ud t,
	\end{equation}
	where $f^{ij}$ is the (assumed constant) fraction of people from $i$ currently in region $j$ at any moment in time.
	
	Equation \eqref{def:flux} is obtained by summing the number of infections in population $n$ due to population $k$ across each region. The summand represents these infections occurring in region $\ell$. This can be seen as the term $f^{n\ell} S^n(t)$ is the proportion of population $n$ that are susceptible and within region $\ell$. Of this population, there will be $\beta f^{n\ell} S^n(t)$ infections per infectious individual per day. Since the number of infectious from $k$ that are in $\ell$ is $f^{k\ell}I^k(t)P^k$, we have that $(\beta f^{n\ell} S^n(t))(f^{k\ell}I^k(t)P^k) \ud t$ is the number of new infections in population $n$ due to population $k$ occurring in the region $\ell$ within the time interval $(t,t+\ud t)$.
	
	Let us assume for now that $n \neq k$. Since we assume that $1 > f^{nn} \gg f^{nk}$, the terms in \eqref{def:flux} besides the cases where $\ell = n$ or $\ell = k$ are negligible. Removing the negligible terms and dividing by the population $P^n$, we see that the influx to the scaled exposed population $E^n$ due to transmission from population $k$ over the interval $(t,t+\ud t)$ is
	\begin{equation}
	\frac{P^k}{P^n} \beta (f^{nn} f^{kn} + f^{nk}f^{kk})S^n(t)I^k(t)\ud t,
	\end{equation}
	which is exactly the influx represented by the model of  $\beta^{nk}S^n(t)I^k(t)\ud t$. This verifies the form of $\beta^{nk}$ for $n \neq k$ in \eqref{def:Beta}. Similarly if we take $n = k$ in \eqref{def:flux} and ignore all terms other than $\ell = n ( = k)$, then we derive the formula for $\beta^{nn}$ in \eqref{def:Beta}.

\section{Detailed discussions on the enhanced DFP algorithm}

\subsection{Implementation details}\label{app:implementation}
In the NY-NJ-PA case study, we choose feedforward architectures for both $V$-networks and $\alpha$-networks. Both have three hidden layers with a width of 40 neurons. The activation function in each hidden layer is $\tanh(x)$. We do not apply activation function to the output layer of $V$-networks, and choose sigmoid function $\rho_s(x) = \frac{1}{1+e^{-x}}$ for the $\alpha$-networks. Other hyperparameters are summarized in the table below.
	\begin{table}[htp!]
		\centering
		\begin{tabular}{ccccccc}
			\hline
			hyperparameter & $lr$ & $M$   & $N_{\text{SGD\_per\_stage}}$ & $N_{\text{batch}}$ & $N_T$ & $\tau$  \\ \hline
			value           & 5e-4  & 250 & 100       & 256       &40      & $1e^{-3}$/180   \\ \hline
		\end{tabular}
		\caption{Hyperparameters in the case study: $lr$ denotes the learning rate in stochastic gradient descent method, $M$ is the total stages of fictitious play, $N_{\text{SGD\_per\_stage}}$ is the number of stochastic gradient descent done in each minimization problem \eqref{eq:cost_discrete}, $N_\text{batch}$ is batch size in each stochastic gradient descent, $N_T$ is the discretization steps on $[0,T]$, and $\tau$ is the weight of the control part in the loss function \eqref{eq:cost_discrete}.}\label{hyper-parameters}
	\end{table}

\subsection{Discussion on the choice of $M$ and $N_{\text{SGD\_per\_stage}}$}\label{app:morenumerics}
\begin{figure}[!ht]
		\centering
		\includegraphics[width =\textwidth, trim = {3em 6em 6em 8em}, clip, keepaspectratio = True]{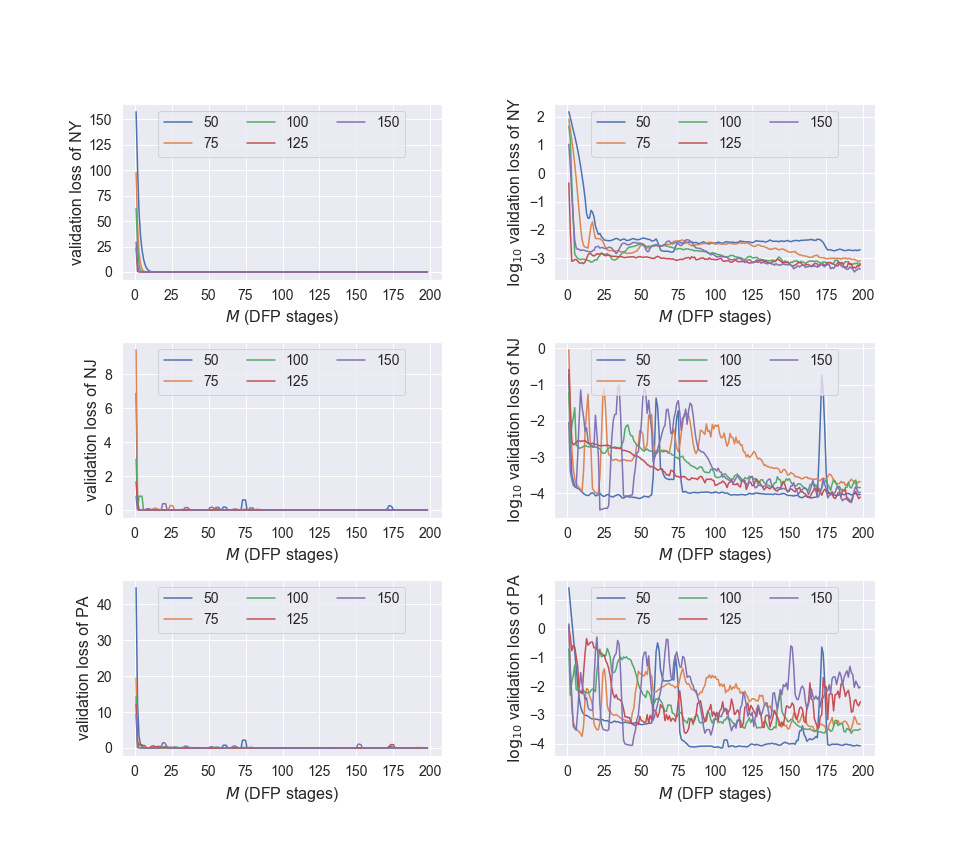}
		\vspace{-2.4em}
		\caption{Loss curves of each state. Left: validation losses versus rounds  $M$ of the enhanced deep fictitious play; Right: $\log_{10}$ validation loss versus rounds $M$ of the enhanced deep fictitious play. The loss curves with respect to $N_{\text{SGD\_per\_stage}}=50,75,100,125,150$ are depicted in blue, orange, green, purple and red. A smoothed moving average with window
			size 3 is applied in the final plots.}
		\label{fig:loss_curve}
	\end{figure}

We provide further experiments here on various choices of $M$ and $N_{\text{SGD\_per\_stage}}$. In Figure~\ref{fig:loss_curve}, we plot both validation loss and log loss against $M$ for three states, which are produced by evaluating the NNs using unseen data after each fictitious play stage. In each panel, loss curves associated with different number of SGDs per stage are presented in different colors (blue = 50, orange = 75, green = 100, red = 125, purple = 150). 
	
	The numerical results show that the validation losses for all states decrease as the number of DFP stages $M$ increases. Moreover, it shows that different $N_{\text{SGD\_per\_stage}}$ generates loss curves with similar patterns. Smaller $N_{\text{SGD\_per\_stage}}$ is more stable on the validation loss of PA. This result is consistent with \cite{HaHu:19} and \cite{seale2006solving}, which convey that it is unnecessary to solve the problem extremely accurate in each stage and that a moderate number of $N_{\text{SGD\_per\_stage}}$ is sufficient. As a result, we choose $N_{\text{SGD\_per\_stage}} = 100$ in our case study.
	
	\subsection{Stability over different experiments}\label{app:NE}
\begin{figure}[!ht]
			\centering
			\includegraphics[width =\textwidth]{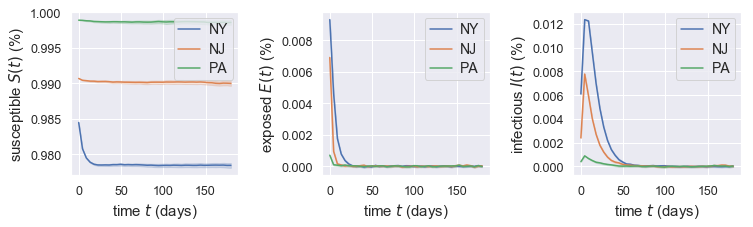}
			\vspace{-2.4em}
			\caption{With the parameter combination $\theta=0.99$, $a=100$, the algorithm identifies one Nash equilibrium for the NY-NJ-PA case study.}
			\label{fig:single-equi}
		\end{figure}

Here we present experiments to investigate the Nash equilibrium of the model with different combinations of parameters. For each combination of parameters, we use the same hyper-parameters and repeat the experiments several times. We run the algorithm for a certain computational budget, and then filter out the results with a fluctuating loss near the stopping and check the converged equilibrium. In the first combination of parameters, we take $\theta=0.99$ and $a=100$, corresponding to the case that a governor weighs the deaths much more than the economic loss and tries to avoid it, and the residents have a strong willingness to follow the governor's policies. Intuitively, the pandemic is possible to get well-controlled. Our numerical experiments confirm this intuition: all converging trails lead to the same Nash equilibrium.  
		A representative plot of $X(t)=(S(t),E(t),I(t))$ is shown in Figure \ref{fig:single-equi}. 
		
		In the second batch of experiments, we take $\theta=0.9$ and $a=50$, corresponding to the case that a governor pays less attention to the number of death and the residents are less willing to follow the policies compared to the first batch of experiments. The change leads to the possibility of multiple Nash equilibrium and the pandemic being out of control. In this case, with different NNs' initialization, the algorithm identifies two Nash equilibria: 
		$75\%$ of the experiments converge to the Nash equilibrium that the pandemic gets controlled and $25\%$ of the experiments converge to the other Nash equilibrium where the pandemic gets out of control.
		\begin{figure}[!ht]
			\centering
			\includegraphics[width =\textwidth]{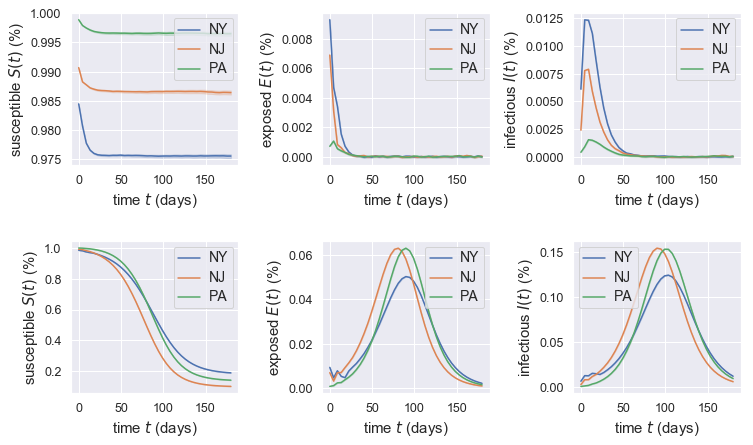}
			\vspace{-2.4em}
			\caption{With the parameter combination $\theta=0.9$, $a=50$, the algorithm identifies two possible Nash equilibria: an under control one (topic panels, with 75\% of the experiments) and on out-of-control one (bottom panels, with 25\% of the experiments).} 
			\label{fig:multi-equi}
		\end{figure}
		
		In conclusion, it is possible to have multiple Nash equilibria depending on the parameter chosen in our stochastic multi-region SEIR model. There is usually a single Nash equilibrium for parameters chosen at extreme values, while for the parameters selected in the middle range, there could exist multiple Nash equilibria. When multiple equilibria exist, we conjecture that the possibility to reach a particular one depends on where we start the fictitious play (the initialization of the NNs' parameters).

\end{document}